\documentclass[leqno,a4paper,11pt]{article}
\usepackage{amsmath}
\usepackage{amssymb}
\usepackage{amsthm}
\usepackage{amscd}

\begin{document}

\newtheorem{thm}{\sc{Theorem}}[section]
\newtheorem{prop}[thm]{\sc{Proposition}}
\newtheorem{lem}[thm]{\sc{Lemma}}
\newtheorem{cor}[thm]{\sc{Corollary}}
\newtheorem{defn}[thm]{\sc{Definition}}
\newtheorem{rem}[thm]{\it{Remark}}
\newtheorem{eg}[thm]{\it{Example}}
\renewcommand{\proofname}{\sc{Proof}}

\numberwithin{equation}{section}

\renewcommand{\theenumi}{\roman{enumi}}
\renewcommand{\labelenumi}{$\mathrm{(\theenumi)}$}

\renewcommand{\refname}{\begin{center} References \end{center}}

\renewcommand{\abstractname}{}

\renewcommand{\thefootnote}{}

\renewcommand{\figurename}{\bf {Fig.}}

\setcounter{section}{-1}

\title{The $b$-functions for prehomogeneous vector spaces of commutative parabolic type and generalized universal Verma modules}

\author{Atsushi KAMITA\\
\\
\textit{Department of Mathematics}\\
\textit{Osaka City University}\\
\textit{Osaka, 558-8585 Japan}\\
E-mail: kamita@sci.osaka-cu.ac.jp
}

\date{}

\maketitle

\begin{abstract}
\hspace{-7mm} \textbf{\sc{Abstract.}}
We shall give a uniform expression and a uniform calculation for the $b$-functions of prehomogeneous vector spaces of commutative parabolic type, which were previously calculated by case-by-case analysis.
Our method is a generalization of Kashiwara's approach using the universal Verma modules.
We shall also give a new proof for the criterion of the irreducibility of the generalized Verma module in terms of $b$-functions due to Suga \cite{Suga}, Gyoja \cite{Gyoja}, Wachi \cite{Wachi}.
\end{abstract}

\section{Introduction}

In this paper we deal with the $b$-functions of the invariants on the flag manifolds $G/P$.
In the case where $P$ is a Borel subgroup, Kashiwara \cite{Kashiwara} determined the $b$-functions by using the universal Verma modules.
Our method is a generalization of Kashiwara's approach.

Let $\mathfrak{g}$ be a simple Lie algebra over the complex number field
$\mathbb{C}$, and let $G$ be a connected simply-connected simple algebraic group with Lie algebra $\mathfrak{g}$.
Fix a parabolic subalgebra $\mathfrak{p}$ of $\mathfrak{g}$.
We denote the reductive part of $\mathfrak{p}$ and the nilpotent part of $\mathfrak{p}$ by $\mathfrak{l}$ and $\mathfrak{n}$ respectively.
Let $L$ be the subgroup of $G$ corresponding to $\mathfrak{l}$.
Let $R$ be the symmetric algebra of the commutative Lie algebra $\mathfrak{p}/[\mathfrak{p}, \mathfrak{p}]$.
For a Lie algebra $\mathfrak{a}$ we set $U_{R}(\mathfrak{a}) = R \otimes_{\mathbb{C}} U(\mathfrak{a})$ where $U(\mathfrak{a})$ denotes the enveloping algebra of $\mathfrak{a}$.
The canonical map $\mathfrak{p} \to R$ induces a one-dimensional $U_R(\mathfrak{p})$-module 
$R_c$.
Let $\mathbb{C}_{\mu}$ be the one-dimensional $\mathfrak{p}$-module with weight $\mu$.
Set $R_{c+\mu} = R_{c} \otimes_{\mathbb{C}} \mathbb{C}_{\mu}$.
Then $R_{c+\mu}$ is a one-dimensional $U_{R}(\mathfrak{p})$-module.

For a character $\mu$ of $\mathfrak{p}$ we regard $\mu$ as a weight of $\mathfrak{g}$, and let $V(\mu)$ be the irreducible $\mathfrak{g}$-module with highest weight $\mu$.
We assume that the weight $\mu$ of $\mathfrak{g}$ is dominant integral.
We define a $U_{R}(\mathfrak{g})$-module homomorphism 
\[
\iota : U_R(\mathfrak{g}) \otimes_{U_R(\mathfrak{p})} R_{c+\mu} \to U_R(\mathfrak{g}) \otimes_{U_R(\mathfrak{p})} (R_{c} \otimes_{\mathbb{C}} V(\mu))
\]
by $\iota(1 \otimes 1) = 1 \otimes 1 \otimes v_{\mu}$, where $v_{\mu}$ is the highest weight vector of $V(\mu)$.
For a $U_{R}(\mathfrak{g})$-module homomorphism $\psi$ from $U_R(\mathfrak{g}) \otimes_{U_R(\mathfrak{p})} (R_{c} \otimes_{\mathbb{C}} V(\mu))$ to $U_R(\mathfrak{g}) \otimes_{U_R(\mathfrak{p})} R_{c+\mu}$ the composite $\psi \iota$ is the multiplication on $U_R(\mathfrak{g}) \otimes_{U_R(\mathfrak{p})} R_{c+\mu}$ by an element $\xi$ of $R$:
\begin{align}\label{eq-intro_cd}
\begin{CD}
U_R(\mathfrak{g}) \otimes_{U_R(\mathfrak{p})} R_{c+\mu} @= U_R(\mathfrak{g}) \otimes_{U_R(\mathfrak{p})} R_{c+\mu}\\
@V{\iota}VV @VV{\xi \, \mathrm{id}}V\\
U_R(\mathfrak{g}) \otimes_{U_R(\mathfrak{p})} (R_{c} \otimes_{\mathbb{C}} V(\mu)) @>>{\psi}> U_R(\mathfrak{g}) \otimes_{U_R(\mathfrak{p})} R_{c+\mu}.
\end{CD}
\end{align}
The set $\Xi_{\mu}$ consisting of all $\xi \in R$ induced by homomorphisms from $U_R(\mathfrak{g}) \otimes_{U_R(\mathfrak{p})} (R_{c} \otimes_{\mathbb{C}} V(\mu))$ to $U_R(\mathfrak{g}) \otimes_{U_R(\mathfrak{p})} R_{c+\mu}$ as above is an ideal of $R$.
We can construct a particular element $\psi_{\mu} \in \mathrm{Hom}_{U_R(\mathfrak{g})}\left(U_R(\mathfrak{g}) \otimes_{U_R(\mathfrak{p})} (R_{c} \otimes_{\mathbb{C}} V(\mu)), U_R(\mathfrak{g}) \otimes_{U_R(\mathfrak{p})} R_{c+\mu} \right)$ by considering the irreducible decomposition of $V(\mu)$ as a $\mathfrak{p}$-module (see Section \ref{sc-UGV} below), however, the corresponding 
$\xi_{\mu} \in \Xi_{\mu}$ is not a generator of $\Xi_{\mu}$ in general.
Note that Kashiwara \cite{Kashiwara} gave the generator of $\Xi_{\mu}$ when $P$ is a Borel subgroup.

Let $\psi \in \mathrm{Hom}_{U_R(\mathfrak{g})}(U_R(\mathfrak{g}) \otimes_{U_R(\mathfrak{p})} (R_{c} \otimes_{\mathbb{C}} V(\mu)), U_R(\mathfrak{g}) \otimes_{U_R(\mathfrak{p})} R_{c+\mu})$ and let $\xi \in \Xi_{\mu}$ be the corresponding element. Then as in Kashiwara \cite{Kashiwara} we can define a differential operator $P(\psi)$ on $G$ satisfying
\[
P(\psi) f^{\lambda + \mu} = \xi(\lambda) f^{\lambda}
\]
for any character $\lambda$ of $\mathfrak{p}$ which can be regarded as a dominant integral weight of $\mathfrak{g}$.
Here, $f^{\lambda}$ denotes the invariant on $G$ corresponding to $\lambda$ (see Section \ref{sc-semi-inv} below) and $\xi$ is regarded as a function
on $\mathrm{Hom}(\mathfrak{p}, \mathbb{C})$.

In the rest of Introduction we assume that the nilpotent radical $\mathfrak{n}$ of $\mathfrak{p}$ is commutative. Then the pair $(L, \mathfrak{n})$ is a prehomogeneous vector space via the adjoint action of $L$.
In this case there exists exactly one simple root $\alpha_0$ such that the root space $\mathfrak{g}_{\alpha_0}$ is in $\mathfrak{n}$.
We denote the fundamental weight corresponding to $\alpha_0$ by $\varpi_0$.

We define an element $\xi_0 \in R$ by
\[
\xi_0(\lambda) = \prod_{\eta \in Wt(\varpi_0) \setminus \{\varpi_0\}} 
\big((\lambda + \rho + \varpi_0, \lambda + \rho + \varpi_0) - (\lambda + \rho + \eta, \lambda + \rho + \eta)\big) \quad (\lambda \in \mathbb{C} \varpi_0),
\]
where $Wt(\varpi_0)$ is the set of the highest weights of irreducible $\mathfrak{l}$-submodules of $V(\varpi_0)$, and $\rho$ is the half sum of positive roots of $\mathfrak{g}$.

\begin{thm}
We have $\xi_{0} = \xi_{\varpi_0}$, and the ideal $\Xi_{\varpi_0}$ of $R$ is generated by $\xi_{0}$.
\end{thm}
We denote by $\psi_0$ the homomorphism satisfying $\psi_0 \iota = \xi_0 \mathrm{id}$.

Let $\mathfrak{n}^-$ be the nilpotent part of the parabolic subalgebra of $\mathfrak{g}$ opposite to $\mathfrak{p}$.
We can define a constant coefficient differential operator $P'(\psi_0)$ on $\mathfrak{n}^- \simeq \exp (\mathfrak{n}^-)$ by
\[
(P(\psi_0) f)|_{\exp (\mathfrak{n}^-)} = P'(\psi_0) (f|_{\exp (\mathfrak{n}^-)}).
\]

\begin{thm}
If the prehomogeneous vector space $(L, \mathfrak{n})$ is regular, then $P'(\psi_0)$ corresponds to the unique irreducible relative invariant of $(L, \mathfrak{n})$, and $b(s) = \xi_0(s \varpi_0)$ is its $b$-function.
\end{thm}

Moreover, using the commutative diagram \eqref{eq-intro_cd} for $\xi_0$ and $\psi_0$
we give a new proof of the following criterion of
the irreducibility of the generalized Verma module due to Suga \cite{Suga}, Gyoja \cite{Gyoja}, Wachi \cite{Wachi}:
\[
U(\mathfrak{g}) \otimes_{U(\mathfrak{p})} \mathbb{C}_{s_0 \varpi_{0}}
\textrm{ is irreducible } \iff \xi_0((s_0 - m)\varpi_{0}) \neq 0 \textrm{ for any positive integer } m.
\]

The author expresses the gratitude to 
Professor T.\ Tanisaki for his valuable advice.

 \section{Prehomogeneous Vector Spaces}
In this section we recall some basic facts on prehomogeneous vector spaces
(see Sato and Kimura \cite{Sato}).

\begin{defn}
\begin{enumerate}
	\item For a connected algebraic group $G$ over the complex number field $\mathbb{C}$ and a finite dimensional $G$-module $V$, the pair $(G, V)$ is called a prehomogeneous vector space if there exists a Zariski open orbit in $V$.
	\item We denote the ring of polynomial functions on $V$ by $\mathbb{C}[V]$.
A nonzero element $f \in \mathbb{C}[V]$ is called a relative invariant
of a prehomogeneous vector space $(G, V)$
if there exists a character $\chi$ of $G$ such that
$f(g v) = \chi(g) f(v)$ for any $g \in G$ and $v \in V$.
	\item A prehomogeneous vector space is called regular if there exists a relative invariant $f$ such that the Hessian $H_f =
\det (\partial^{2} f / \partial x_i \partial x_j)$ is not identically zero,
where $\{ x_i \}$ is a coordinate system of $V$.
\end{enumerate}
\end{defn}
We call algebraically independent relative invariants $f_{1}, f_2, \ldots, f_l$
basic relative invariants if for any relative invariant $f$ there exist $c \in \mathbb{C}$ and $m_i \in \mathbb{Z}$ such that $f = c f_1^{m_1} \dotsm f_{l}^{m_l}$.

Assume that $(G, V)$ is a prehomogeneous vector space such that $G$ is reductive.
Then the dual space $V^{*}$ of $V$ is also a prehomogeneous vector space
by $\langle g v^*, v \rangle = \langle v^*, g^{-1} v \rangle$,
where $\langle \ , \ \rangle$ is the natural pairing of $V^{*}$ and $V$.
If $f \in \mathbb{C}[V]$ is a relative invariant of $(G, V)$ with character $\chi$, then there exists a relative invariant $f^{*}$ of $(G, V^{*})$ with character $\chi^{-1}$.
For $h \in \mathbb{C}[V^*]$ we define a constant coefficient differential operator $h(\partial)$ by
\[
h(\partial) \exp \langle v^{*}, v \rangle = h(v^{*}) \exp \langle v^{*} v \rangle,
\]
where $v \in V$ and $v^* \in V^*$.
Then there exists a polynomial $b(s) \in \mathbb{C}[s]$ such that
\[
f^{*}(\partial) f^{s + 1} = b(s) f^s.
\]
This polynomial is called the $b$-function of $f$.
It is known that $\deg b = \deg f$ (see \cite{Sato0}).

\section{Generalized Universal Verma Modules}\label{sc-UGV}
Let $\mathfrak{g}$ be a simple Lie algebra over $\mathbb{C}$ with Cartan subalgebra $\mathfrak{h}$.
Let $\Delta \subset \mathfrak{h}^{*}$ be the root system and
$W \subset \mathrm{GL}(\mathfrak{h})$ the Weyl group.
For $\alpha \in \Delta$ we denote the corresponding root space by
$\mathfrak{g}_{\alpha}.$
We denote the set of positive roots by $\Delta^{+}$
and the set of simple roots by $\{\alpha_{i}\}_{i \in I_0}$,
where $I_0$ is an index set.
Let $\rho$ be the half sum of positive roots of $\mathfrak{g}$.
We set
\[
\mathfrak{n}^{\pm} = \bigoplus_{\alpha \in \Delta^{+}} \mathfrak{g}_{\pm \alpha},
\hspace{5mm}
\mathfrak{b}^{\pm} = \mathfrak{h} \oplus \mathfrak{n}^{\pm}.
\]
For $i \in I_0$ let $h_i \in \mathfrak{h}$ be the simple coroot and
$\varpi_i \in \mathfrak{h}^{*}$ the fundamental weight corresponding to $i$.
We denote the longest element of $W$ by $w_0$.
Let $(\ , \ )$ be the $W$-invariant nondegenerate symmetric bilinear form
on $\mathfrak{h}^{*}$.
We denote the irreducible $\mathfrak{g}$-module with highest weight $\mu \in \sum_{i \in I_0} \mathbb{Z}_{\ge 0} \varpi_{i}$ 
by $V(\mu)$ and its highest weight vector
by $v_{\mu}$. 
For a Lie algebra $\mathfrak{a}$ we denote the enveloping algebra of
$\mathfrak{a}$ by $U(\mathfrak{a})$.

For a subset $I \subset I_0$ we set
\begin{align*}
&\Delta_I = \Delta \cap \sum_{i \in I} \mathbb{Z} \alpha_i,&
&\mathfrak{l}_I = \mathfrak{h} \oplus (\bigoplus_{\alpha \in \Delta_I} \mathfrak{g}_{\alpha}),&\\
&\mathfrak{n}_{I}^{\pm} = \bigoplus_{\alpha \in \Delta^{+} \setminus \Delta_I} \mathfrak{g}_{\pm \alpha},&
&\mathfrak{p}_I^{\pm} = \mathfrak{l}_I \oplus \mathfrak{n}_I^{\pm},&\\
&\mathfrak{h}_I = \mathfrak{h} / \sum_{i \in I} \mathbb{C} h_i,&
&\mathfrak{h}_I^{*} = \sum_{i \in I_{0} \setminus I} \mathbb{C} \varpi_i.&
\end{align*}
Let $W_I$ be the subgroup of $W$ generated by the simple reflections corresponding
to $i \in I$.
We denote the longest element of $W_I$ by $w_I$.
Let $\mathfrak{h}_{I, +}^{*}$ be the set of dominant integral weights in
$\mathfrak{h}_I^{*}$.
For $\mu \in \mathfrak{h}_I^*$ we define a one-dimensional $U(\mathfrak{p}_I^+)$-module
$\mathbb{C}_{I, \mu}$ by
\[
\mathbb{C}_{I, \mu} = U(\mathfrak{p}_I^+) \Big/ \big(U(\mathfrak{p}_I^{+}) \mathfrak{n}^{+} + \sum_{h \in \mathfrak{h}} U(\mathfrak{p}_I^{+}) (h - \mu(h)) + U(\mathfrak{p}_I^{+}) (\mathfrak{n}^{-} \cap \mathfrak{l}_I)\big).
\]
We denote the canonical generator of $\mathbb{C}_{I, \mu}$ by $1_{\mu}$.
Set $M_I(\mu) = U(\mathfrak{g}) \otimes_{U(\mathfrak{p}_I^+)} \mathbb{C}_{I, \mu}$.
We denote the irreducible $\mathfrak{p}_I^{+}$-module with highest weight $\mu \in \sum_{i \in I} \mathbb{Z}_{\ge 0} \varpi_{i} + \sum_{j \notin I} \mathbb{Z} \varpi_{j}$ by $W(\mu)$.

Let $G$ be a connected simply-connected simple algebraic group with Lie algebra $\mathfrak{g}$.
We denote the subgroups of $G$ corresponding to $\mathfrak{h}, \mathfrak{b}^{\pm}, \mathfrak{l}_I, \mathfrak{n}_I^{\pm}$ by $T, B^{\pm}, L_I, N_I^{\pm}$ respectively.

Let $R_I$ be the symmetric algebra of $\mathfrak{h}_I$,
and define a linear map $c : \mathfrak{h} \to R_I$ as the composite of
the natural projection from $\mathfrak{h}$ to $\mathfrak{h}_I$
and the natural injection from $\mathfrak{h}_I$ to $R_I$.
Set $U_{R_I}(\mathfrak{a}) = R_I \otimes_{\mathbb{C}} U(\mathfrak{a})$ for a Lie algebra $\mathfrak{a}$.

We set for $\mu \in \mathfrak{h}_I^{*}$
\begin{align*}
R_{I, c + \mu} = U_{R_I} (\mathfrak{p}_I^{+}) \Big/ \big(U_{R_I}(\mathfrak{p}_I^{+}) \mathfrak{n}^{+} + \sum_{h \in \mathfrak{h}} U_{R_I}(\mathfrak{p}_I^{+}) (h - c(h) - \mu(h)) + U_{R_I}(\mathfrak{p}_I^{+}) (\mathfrak{n}^{-} \cap \mathfrak{l}_I)\big).
\end{align*}
We denote the canonical generator of $R_{I, c + \mu}$ by $1_{c + \mu}$.
\begin{defn}\label{def-universal Verma}
For $\mu \in \mathfrak{h}_I^{*}$ we call $M_{R_I}(c + \mu) = 
U_{R_I}(\mathfrak{g}) \otimes_{U_{R_I}(\mathfrak{p}_I^+)} R_{I, c + \mu}$
a generalized universal Verma module.
\end{defn}
Note that $M_{R_{\emptyset}}(c)$ is the universal Verma module in
Kashiwara \cite{Kashiwara}.
For $\lambda \in \mathfrak{h}_I^*$ we regard $\mathbb{C}$
as an $R_I$-module by $c(h_i) 1 = \lambda(h_i)$.
Then we have 
\[
\mathbb{C} \otimes_{R_I} M_{R_I}(c + \mu) = M_I(\lambda + \mu).
\]

The next lemma is obvious.
\begin{lem}
$\mathrm{End}_{U_{R_I}(\mathfrak{g})}(M_{R_I}(c + \mu)) = R_I$.
\end{lem}

For $\mu \in \mathfrak{h}_{I}^{*}$ we define a $U_{R_I}(\mathfrak{g})$-module homomorphism
\[
\iota_{\mu} : M_{R_I}(c + \mu) \longrightarrow U_{R_I}(\mathfrak{g}) \otimes_{U_{R_I}(\mathfrak{p}_I^+)}(R_{I,c} \otimes_{\mathbb{C}} V({\mu}))
\]
by $\iota_{\mu}(1 \otimes 1_{c + \mu}) = 1 \otimes 1_{c} \otimes v_{\mu}$.
We denote by $\Xi_{\mu}$ the ideal of $R_I$ consisting of $\xi$ such that there exists a homomorphism $\psi \in \mathrm{Hom}_{U_{R_I}(\mathfrak{g})}(U_{R_I}(\mathfrak{g}) \otimes_{U_{R_I}(\mathfrak{p}_I^+)}(R_{I,c} \otimes_{\mathbb{C}} V({\mu})), M_{R_I}(c + \mu))$ satisfying $\psi \iota_{\mu} = \xi \, \mathrm{id}$.
Let us give a particular element $\xi_{\mu}$ of $\Xi_{\mu}$ for  $\mu \in \mathfrak{h}_{I, +}^*$.

\begin{lem}\label{lem-Ext}
For $\mu_1, \mu_2 \in \sum_{i \in I} \mathbb{Z}_{\ge 0} \varpi_{i} + \sum_{j \notin I} \mathbb{Z} \varpi_{j}$ we define a function $p_{\mu_1, \mu_2}$
on $\mathfrak{h}_I^*$ by
\[
p_{\mu_1, \mu_2}(\lambda) = (\lambda + \rho + \mu_1, \lambda + \rho + \mu_1) -
(\lambda + \rho + \mu_2, \lambda + \rho + \mu_2),
\]
which is regarded as an element of $R_I$.
Then we have
\begin{align*}
p_{\mu_1, \mu_2} \ \mathrm{Ext}_{U_{R_I}(\mathfrak{g})}^1
\big(U_{R_I}(\mathfrak{g}) \otimes_{U_{R_I}(\mathfrak{p}_I^{+})} (R_{I,c} \otimes_{\mathbb{C}} W(\mu_1)), \ 
U_{R_I}(\mathfrak{g}) \otimes_{U_{R_I}(\mathfrak{p}_I^{+})} (R_{I,c} \otimes_{\mathbb{C}} W(\mu_2))\big) \\ 
= 0.
\end{align*}
\end{lem}

\begin{proof}
The action of the Casimir element of $U(\mathfrak{g})$ on $U_{R_I}(\mathfrak{g}) \otimes_{U_{R_I}(\mathfrak{p}_I^{+})} (R_{I,c} \otimes_{\mathbb{C}} W(\mu))$ is given by the multiplication by $p_{\mu} \in R_I$, where $p_{\mu}(\lambda) = (\lambda + \rho + \mu, \lambda + \rho + \mu) - (\rho, \rho)$ for $\lambda \in \mathfrak{h}_I^*$.
Using this action, we can easily check that $p_{\mu_1, \mu_2} = p_{\mu_1} - p_{\mu_2}$ is an annihilator.
\end{proof}

\begin{lem}\label{lem-decomp of p-module}
For any $\mu \in \mathfrak{h}_{I, +}^{*}$
there exist $\mathfrak{p}_I^+$-submodules $F_1, F_2, \ldots, F_r$ of $V(\mu)$
and weights $\eta_1, \eta_2, \ldots, \eta_{r-1} \in \sum_{i \in I} \mathbb{Z}_{\ge 0} \varpi_i + \sum_{i \in I_0 \setminus I} \mathbb{Z} \varpi_{i}$
satisfying the following conditions:
\begin{enumerate}
	\item $\mathbb{C} v_{\mu} = F_1 \subsetneq F_2 \subsetneq \cdots
	\subsetneq F_r = V(\mu)$.
	\item $F_{i+1} / F_i \simeq W(\eta_i)^{\oplus N_i}$ for some positive
	integer $N_i$.
	\item $\eta_i \neq \eta_j$ for $i \neq j$.
\end{enumerate}
\end{lem}

\begin{proof}
For a non-negative integer $m$ we set
\[
P(m) = \{ \lambda \in \mathfrak{h}^* \ | \  \mu - \lambda = \sum_{i \in I_0} m_i \alpha_i \textrm{ and } \sum_{i \notin I} m_i = m \} \ \textrm{ and } \ 
V_m = \bigoplus_{\lambda \in P(m)} V(\mu)_{\lambda},
\]
where $V(\mu)_{\lambda}$ is the weight space of $V(\mu)$ with weight $\lambda$.
Then $V_m$ is an $\mathfrak{l}_I$-module,
and we have the irreducible decomposition
\[
V_m = \tilde{W}(\eta_{m,1})^{\oplus N_{m,1}} \oplus \cdots \oplus \tilde{W}(\eta_{m,t_m})^{\oplus N_{m,t_m}}
\]
where $\tilde{W}(\eta)$ is the irreducible $\mathfrak{l}_I$-module with highest
weight $\eta$, and $\eta_{m,i} \neq \eta_{m, j}$ for $i \neq j$.
For $1 \le i \le t_m$ we define a $\mathfrak{p}_I^{+}$-submodule $F_{m, i}$
of $V(\mu)$ by
\[
F_{m, i} = V_0 \oplus \cdots \oplus V_{m-1} \oplus \tilde{W}(\eta_{m,1})^{\oplus N_{m,1}} \oplus \cdots \oplus \tilde{W}(\eta_{m,i})^{\oplus N_{m,i}}.
\]
Then we have the sequence
\begin{align*}
\mathbb{C}v_{\mu} = F_{0,1} \subsetneq \cdots \subsetneq F_{m-1, t_{m-1}} \subsetneq F_{m,1} \subsetneq F_{m,2} \subsetneq \cdots \subsetneq F_{m, t_m}\subsetneq \cdots \subsetneq F_{r, t_r} = V(\mu).
\end{align*}
It is clear that the above sequence satisfies the conditions (ii) and (iii).
\end{proof}

For $\mu \in \mathfrak{h}_{I, +}^{*}$, we fix the sequence
$\{F_1, F_2, \ldots, F_r\}$ of $\mathfrak{p}_I^{+}$-submodules of $V(\mu)$
satisfying the conditions of Lemma \ref{lem-decomp of p-module},
and set $\xi_{\mu} = \prod_{i=1}^{r-1} p_{\mu, \eta_i} \in R_I$.

\begin{thm}\label{thm-commutative diag}
For $\mu \in \mathfrak{h}_{I,+}^*$ we have $\xi_{\mu} \in \Xi_{\mu}$.
\end{thm}
\begin{proof}
It is clear that
$U_{R_I}(\mathfrak{g}) \otimes_{U_{R_I}(\mathfrak{p}_I^+)}(R_{I,c} \otimes_{\mathbb{C}} F_1) \simeq M_{R_I}(c + \mu)$.
Let $\iota_j$ be the canonical injection from $U_{R_I}(\mathfrak{g}) \otimes_{U_{R_I}(\mathfrak{p}_I^+)}(R_{I,c} \otimes_{\mathbb{C}} F_{j})$ into
$U_{R_I}(\mathfrak{g}) \otimes_{U_{R_I}(\mathfrak{p}_I^+)}(R_{I,c} \otimes_{\mathbb{C}} F_{j+1})$.
We show that there exists a commutative diagram
\begin{equation}
\begin{CD}\label{eq-commutative diag large}
U_{R_I}(\mathfrak{g}) \otimes_{U_{R_I}(\mathfrak{p}_I^+)}(R_{I,c} \otimes_{\mathbb{C}} F_1) @=
 M_{R_I}(c + \mu)\\
@V{\iota_{j-1} \cdots \iota_{1}}VV @VV{\prod_{i = 1}^{j-1}p_{\mu, \eta_i}}V\\
U_{R_I}(\mathfrak{g}) \otimes_{U_{R_I}(\mathfrak{p}_I^+)}(R_{I,c} \otimes_{\mathbb{C}} F_j) @>>{\psi_{j}}> M_{R_I}(c + \mu)
\end{CD}
\end{equation}
by the induction on $j$.
Assume that there exists a commutative diagram \eqref{eq-commutative diag large} for $j \ (\ge 1)$.
From the exact sequence
\begin{align*}
\begin{CD}
0 @>>> U_{R_I}(\mathfrak{g}) \otimes_{U_{R_I}(\mathfrak{p}_I^+)} (R_{I, c} \otimes_{\mathbb{C}} F_j) @>{\iota_j}>> U_{R_I}(\mathfrak{g}) \otimes_{U_{R_I}(\mathfrak{p}_I^+)} (R_{I, c} \otimes_{\mathbb{C}} F_{j+1})\\
\end{CD}\\
\begin{CD}
@>>> U_{R_I}(\mathfrak{g}) \otimes_{U_{R_I}(\mathfrak{p}_I^+)} (R_{I, c} \otimes_{\mathbb{C}} F_{j+1}/F_j) @>>> 0,
\end{CD}
\end{align*}
we have a long exact sequence
\begin{align*}
\begin{CD}
0 @>>> \mathrm{Hom}_{U_{R_I}(\mathfrak{g})}\big(U_{R_I}(\mathfrak{g}) \otimes_{U_{R_I}(\mathfrak{p}_I^+)} (R_{I, c} \otimes_{\mathbb{C}} F_{j+1}/F_j), \ 
M_{R_I}(c + \mu)\big) @. @. \\
@>>> \mathrm{Hom}_{U_{R_I}(\mathfrak{g})}\big(U_{R_I}(\mathfrak{g}) \otimes_{U_{R_I}(\mathfrak{p}_I^+)} (R_{I, c} \otimes_{\mathbb{C}} F_{j+1}), \ 
M_{R_I}(c + \mu)\big) @. @. \\
@>>> \mathrm{Hom}_{U_{R_I}(\mathfrak{g})}\big(U_{R_I}(\mathfrak{g}) \otimes_{U_{R_I}(\mathfrak{p}_I^+)} (R_{I, c} \otimes_{\mathbb{C}} F_j), \ 
M_{R_I}(c + \mu)\big) @. @. \\
@>{\delta}>> \mathrm{Ext}_{U_{R_I}(\mathfrak{g})}^{1}\big(U_{R_I}(\mathfrak{g}) \otimes_{U_{R_I}(\mathfrak{p}_I^+)} (R_{I, c} \otimes_{\mathbb{C}} F_{j+1}/F_j), \ 
M_{R_I}(c + \mu)\big) @>>> \cdots.  
\end{CD}
\end{align*}
Since $F_{j+1}/F_j \simeq W(\eta_j)^{\oplus N_j}$,
we have $\delta(p_{\mu, \eta_j} \psi_j) = p_{\mu, \eta_j} \delta(\psi_j) = 0$
by Lemma \ref{lem-Ext}.
Hence there exists an element $\psi_{j+1} \in \mathrm{Hom}_{U_{R_I}(\mathfrak{g})}\big(U_{R_I}(\mathfrak{g}) \otimes_{U_{R_I}(\mathfrak{p}_I^+)} (R_{I, c} \otimes_{\mathbb{C}} F_{j+1}), \ M_{R_I}(c + \mu)\big)$ such that $\psi_{j+1} \iota_{j} = p_{\mu, \eta_j} \psi_j$.
Hence we have the commutative diagram
\begin{align*}
\begin{CD}
U_{R_I}(\mathfrak{g}) \otimes_{U_{R_I}(\mathfrak{p}_I^+)}(R_{I,c} \otimes_{\mathbb{C}} F_1) @=
 M_{R_I}(c + \mu)\\
@V{\iota_{j-1} \cdots \iota_{1}}VV @VV{\prod_{i = 1}^{j-1}p_{\mu, \eta_i}}V\\
U_{R_I}(\mathfrak{g}) \otimes_{U_{R_I}(\mathfrak{p}_I^+)}(R_{I,c} \otimes_{\mathbb{C}} F_j) @>>{\psi_{j}}> M_{R_I}(c + \mu)\\
@V{\iota_j}VV @VV{p_{\mu, \eta_j}}V\\
U_{R_I}(\mathfrak{g}) \otimes_{U_{R_I}(\mathfrak{p}_I^+)}(R_{I,c} \otimes_{\mathbb{C}} F_{j+1}) @>>{\psi_{j+1}}> M_{R_I}(c + \mu).
\end{CD}
\end{align*}
In particular $\psi_{r} \iota_{\mu} = \xi_{\mu}$.
Therefore $\xi_{\mu} \in \Xi_{\mu}$.
\end{proof}
Let $\psi_{\mu}$ be a homomorphism from $U_{R_I}(\mathfrak{g}) \otimes_{U_{R_I}(\mathfrak{p}_I^+)}(R_{I,c} \otimes_{\mathbb{C}} V(\mu))$ to $M_{R_I}(c + \mu)$
satisfying $\psi_{\mu} \iota_{\mu} = \xi_{\mu}$.
Note that $\psi_{\mu}$ is non-zero since $\xi_{\mu} \neq 0$.

\begin{rem}
\begin{enumerate}
	\item In general $\xi_{\mu}$ is not a generator of the ideal $\Xi_{\mu}$.
For example let $\mathfrak{g}$ be a simple Lie algebra of type $G_2$.
We take the simple roots $\alpha_1$ and $\alpha_2$ such that
$\alpha_1$ is short.
If $I = \{2\}$ and $\mu = \varpi_{1}$, then we have
$\xi_{\mu} = (c(h_1) + 1) (c(h_1) + 2) (c(h_1) + 3) (2c(h_1) + 5)$ up to constant multiple.
But $(c(h_1) + 1) (2 c(h_1) + 5) \in \Xi_{\mu}$.
	\item For $I = \emptyset$ it is shown in Kashiwara \cite{Kashiwara}
that $\Xi_{\mu}$ is generated by
\[
\prod_{\alpha \in \Delta^{+}} \prod_{j = 0}^{\mu(h_{\alpha})-1}(c(h_{\alpha}) + \rho(h_{\alpha}) + j),
\]
where $h_{\alpha}$ is the coroot corresponding to $\alpha$.
\end{enumerate}
\end{rem}

\section{Semi-invariants}\label{sc-semi-inv}
Let $\lambda$ be a dominant integral weight.
We regard the dual space $V(\lambda)^*$ as a left $\mathfrak{g}$-module via
$\langle x v^*, v \rangle$ = $\langle v^*, -x v \rangle$ for $x \in \mathfrak{g}$, $v^* \in V(\lambda)^{*}$ and $v \in V(\lambda)$.
We denote the lowest weight vector of $V(\lambda)^{*}$ by $v_{\lambda}^{*}$.
We normalize $v_{\lambda}^{*}$ by $\langle v_{\lambda}^{*}, v_{\lambda} \rangle = 1$.
\begin{defn}
We define a regular function $f^{\lambda}$ on $G$ by $f^{\lambda}(g) = \langle v_{\lambda}^{*}, g v_{\lambda} \rangle$.
\end{defn}
For $b^{\pm} \in B^{\pm}$ and $g \in G$ we have
\[
f^{\lambda}(b^{-} g b^{+}) = \lambda^{-}(b^{-}) \lambda^+(b^+) f^{\lambda}(g),
\]
where $\lambda^{\pm}$ is the character of $B^{\pm}$ corresponding to $\lambda$.
This function $f^{\lambda}$ is called $B^- \times B^+$-semi-invariant.
Note that $f^{\lambda_1 + \lambda_2} = f^{\lambda_1} f^{\lambda_2}$.

Let $\mu \in \mathfrak{h}_{I, +}^*$.
We take a basis $\{ v_{\mu, j} \}_{0 \le j \le n}$
of $V(\mu)$ consisting of weight vectors such that
$v_{\mu, 0} = v_{\mu}$ is the highest weight vector and
$v_{\mu, n}$ is the lowest.
We denote the dual basis of $V(\mu)^{*}$ by $\{v_{\mu, j}^{*}\}$.
For a $U_{R_I}(\mathfrak{g})$-module homomorphism
\[
\begin{CD}
\psi : U_{R_I}(\mathfrak{g}) \otimes_{U_{R_I}(\mathfrak{p}_{I}^{+})} (R_{I, c} \otimes_{\mathbb{C}} V(\mu)) @>>>  M_{R_I}(c + \mu)
\end{CD}
\]
we define elements $Y'_j \in U_{R_I}(\mathfrak{n}_I^-)$ for $0 \le j \le n$ by
\[
\psi(1 \otimes 1_{c} \otimes v_{\mu, j}) = Y'_{j} \otimes 1_{c + \mu},
\]
and define an element $\xi \in \Xi_{\mu}$ by $\xi = \psi \iota_{\mu}$.
Note that $Y'_{0} = \xi$.
Let $\pi : R_I \to U(\sum_{i \notin I}\mathbb{C}h_i)$ be the algebra isomorphism defined by $\pi(c(h_i)) = h_i - \mu(h_i)$ for $i \notin I$.
Set $\pi(\sum_j a_j \otimes y_j) = \sum_j y_j \pi(a_j)$ for $a_j \in R_I$ and $y_j \in U(\mathfrak{n}_I^-)$.
Clearly we have $y \otimes 1_{c+\mu} = \pi(y) \otimes 1_{c+\mu} \in M_{R_I}(c + \mu)$ \ ($y \in U_{R_I}(\mathfrak{n}_I^-)$).
We set $Y_j = \pi(Y'_j)$.
We define differential operators $P_{\mu}(\psi)$ and $\tilde{P}_{\mu}(\psi)$ on $G$ by
\begin{align*}
&(P_{\mu}(\psi) \varphi)(g) = \sum_{j = 0}^n
\langle g v_{\mu, j}^{*}, v_{\mu, 0} \rangle (R(Y_j)\varphi)(g),\\
&(\tilde{P}_{\mu}(\psi)\varphi)(g) = \sum_{j=0}^n \langle g v_{\mu,j}^*, v_{\mu, n} \rangle (R(Y_j)\varphi)(g),
\end{align*}
where $R(y)$ ($y \in U(\mathfrak{g})$) is the left invariant differential operator induced by the right action of $G$ on itself.
Then we have the following theorem.

\begin{thm}\label{thm-exist b-func}
Let $\mu \in \mathfrak{h}_{I, +}^{*}$ and $\psi \in \mathrm{Hom}_{U_{R_I}(\mathfrak{g})}(U_{R_I}(\mathfrak{g}) \otimes_{U_{R_I}(\mathfrak{p}_{I}^{+})} (R_{I, c} \otimes_{\mathbb{C}} V(\mu)), M_{R_I}(c + \mu))$.
Then we have
\[
P_{\mu}(\psi) f^{\lambda + \mu} = \xi(\lambda) f^{\lambda}
\]
for any $\lambda \in \mathfrak{h}_{I, +}^{*}$.
Here $\xi$ is the element of $\Xi_{\mu}$ defined by $\xi = \psi \iota_{\mu}$.
\end{thm}
We omit the proof since it is almost identical to the one for \cite[Theorem 2.1]{Kashiwara}.

For a dominant integral weight $\lambda$ we define a function $\tilde{f}^{\lambda}$ on $G$ by
\[
\tilde{f}^{\lambda}(g) = \langle v_{w_{0} \lambda}^{*}, g v_{\lambda} \rangle,
\]
where $v_{w_{0} \lambda}^{*}$ is the highest weight vector which is normalized by
$\langle v_{w_{0} \lambda}^{*}, \dot w_{0} v_{\lambda} \rangle = 1$ and 
$\dot w_0 \in N_G(T)$ is a representative element of $w_0 \in W = N_G(T)/T$.
Since $\tilde{f}^{\lambda}(\dot w_{0} g) = f^{\lambda}(g)$, we obtain the following lemma.

\begin{lem}\label{lem-P and tilde}
Let $\lambda, \mu \in \mathfrak{h}_{I,+}^*$. For any $g \in G$
we have $(\tilde{P}_{\mu}(\psi) \tilde{f}^{\lambda})(\dot w_{0} g) = (P_{\mu}(\psi) {f}^{\lambda}) (g)$.
\end{lem}

By Theorem \ref{thm-exist b-func} we have the following corollary.

\begin{cor}\label{cor-b-function}
Let $\mu \in \mathfrak{h}_{I, +}^{*}$. We have
\[
\tilde{P}_{\mu}(\psi) \tilde{f}^{\lambda + \mu} = \xi(\lambda) \tilde{f}^{\lambda}
\]
for any $\lambda \in \mathfrak{h}_{I, +}^{*}$.
Here $\xi$ is the element of $\Xi_{\mu}$ defined by $\xi = \psi \iota_{\mu}$.
\end{cor}

\section{Commutative Parabolic Type}\label{sc-commuta}
In the remainder of this paper we assume that
\[
I = I_{0} \setminus \{ i_0 \}
\]
and that the highest root $\theta$ of $\mathfrak{g}$ is in
$\alpha_{i_0} + \sum_{i \neq i_{0}} \mathbb{Z}_{\ge 0} \alpha_{i}$.
Then it is known that $[\mathfrak{n}_I^{\pm}, \mathfrak{n}_I^{\pm}] = \{0\}$
and the pairs $(L_I, \mathfrak{n}_I^{\pm})$ are prehomogeneous vector spaces via the adjoint action, which are called of commutative parabolic type.
The all pairs $(\mathfrak{g}, i_0)$ of commutative parabolic type are given by
the Dynkin diagrams of Fig. \ref{fig-commutative}.
Here the white vertex corresponds to $i_0$.
\def\Line{-- \hspace{-6mm} -- \hspace{-6mm} --}
{\small \begin{figure}[tbp]
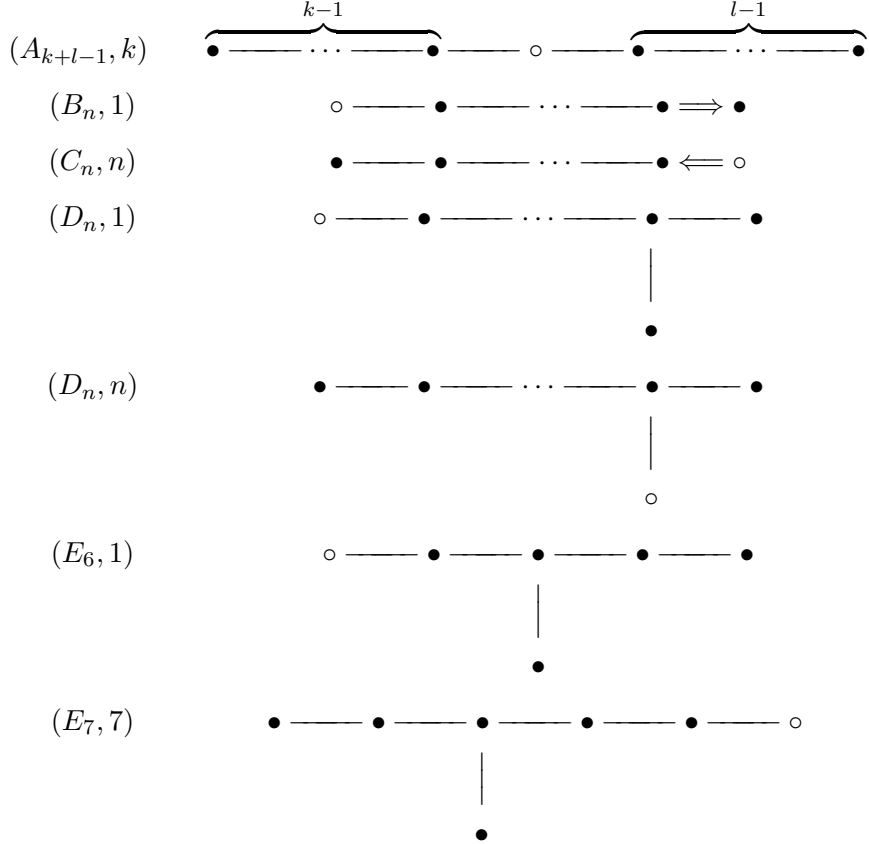

\begin{equation*}
\begin{CD}
(A_{k + l -1}, k) \quad \quad @.
\overbrace{
\bullet \Line \cdots \Line \bullet}^{k-1} \Line \circ
\Line \overbrace{\bullet \Line \cdots \Line \bullet}^{l-1}\\
(B_n, 1) \quad @.
\circ \Line \bullet \Line \cdots \Line \bullet \Longrightarrow \bullet\\
(C_n, n) \quad @.
\bullet \Line \bullet \Line \cdots \Line \bullet \Longleftarrow \circ\\
(D_n, 1) \quad @.
\circ \Line \bullet \Line \cdots \Line \bullet \Line \bullet\\
@.  \hspace{30mm}\Big|\\
@.  \hspace{30mm}\bullet\\
(D_n, n) \quad @.
\bullet \Line \bullet \Line \cdots \Line \bullet \Line \bullet\\
@.  \hspace{30mm}\Big|\\
@.  \hspace{30mm}\circ\\
(E_6, 1) \quad @.
\circ \Line \bullet \Line \bullet \Line \bullet \Line \bullet\\
@.  \Big|\\
@.  \bullet\\
(E_7, 7) \quad @.
\bullet \Line \bullet \Line \bullet \Line \bullet \Line \bullet \Line \circ\\
@.  \hspace{-15mm}\Big|\\
@.  \hspace{-15mm}\bullet\\
\end{CD}
\end{equation*}
\caption{\  Commutative Parabolic Type}
\label{fig-commutative}
\end{figure}}

The pairs $(\mathfrak{g}, i_0)$ such that the corresponding prehomogeneous vector spaces are regular are as follows:
$(A_{2n-1}, n), (B_n, 1), (C_n, n), (D_n, 1), (D_{2n}, 2n)$ and $(E_7,7)$.
Then it is seen that $w_0 \alpha_{i_0} = -\alpha_{i_0}$.

Since $\mathfrak{n}_I^{-}$ is identified with the dual space of
$\mathfrak{n}_I^{+}$ via the Killing form,
the symmetric algebra $S(\mathfrak{n}_I^{-})$ is isomorphic to
$\mathbb{C}[\mathfrak{n}_I^{+}]$.
By the commutativity of $\mathfrak{n}_I^{-}$
we have $S(\mathfrak{n}_I^{-}) = U(\mathfrak{n}_I^{-})$.
Hence $\mathbb{C}[\mathfrak{n}_I^{+}]$ is identified with
$U(\mathfrak{n}_I^{-})$.

Set $\gamma_{1} = \alpha_{i_0}$. For $i \ge 1$ we take the root $\gamma_{i+1}$ as a minimal element in
\[
\Gamma_i = \{ \alpha \in \Delta^{+} \setminus \Delta_{I} \ | \ \alpha + \gamma_{j} \notin \Delta \textrm{ and } \alpha - \gamma_{j} \notin \Delta \cup \{0\} \textrm{ for all } j \le i \}.
\]
Let $r = r(\mathfrak{g})$ be the index such that $\Gamma_{r(\mathfrak{g})-1} \neq \emptyset$ and $\Gamma_{r(\mathfrak{g})} = \emptyset$.
Note that $(\gamma_i, \gamma_j) = 0$ for $i \neq j$.
It is known that all $\gamma_i$ have the same length (see Moore \cite{Moore}).
Set $\mathfrak{h}^- = \sum_{i=1}^r \mathbb{C} h_{\gamma_i}$, where $h_{\gamma_i}$ is the coroot corresponding to $\gamma_i$.
For $1 \le i \le r$ we set $\lambda_{i} = -(\gamma_1 + \cdots + \gamma_i)$.
Then we have the following lemmas.

\begin{lem}[{\rm Moore \cite{Moore}}]\label{lem of gamma-0}
For $\beta \in \Delta^+ \cap \Delta_I$
the restriction $\beta |_{\mathfrak{h}^-}$ is as follows.
\begin{enumerate}
	\item $\beta |_{\mathfrak{h}^-} = 0$. Then $\beta \pm \gamma_i \notin \Delta$ for any $i$.
	\item $\beta |_{\mathfrak{h}^-} = \frac{\gamma_j}{2} |_{\mathfrak{h}^-}$. Then $\beta \pm \gamma_i \notin \Delta$ for any $i \neq j$.
	\item $\beta |_{\mathfrak{h}^-} = \frac{\gamma_j - \gamma_k}{2} |_{\mathfrak{h}^-} \ (j > k)$. Then $\beta \pm \gamma_i \notin \Delta$ for any $i \neq j,k$ and $\beta + \gamma_j \notin \Delta$.
\end{enumerate}
\end{lem}

Set $D = \{\alpha_i \ | \ i \in I\}$.
For a subset $S$ of $\Delta$, $S(\mathfrak{h}^-)$ is defined by
\[
S(\mathfrak{h}^-) = \{ \beta \in \sum_{i = 1}^r \mathbb{Q}\gamma_i \ | \ \beta |_{\mathfrak{h}^-} = \alpha |_{\mathfrak{h}^-} \ \textrm{for $\alpha \in S$} \}.
\]
\begin{lem}[{\rm Moore \cite{Moore}, Wachi \cite{Wachi}}]\label{lem of gamma}
\begin{enumerate}
	\item If $(L_I, \mathfrak{n}_I^{+})$ is regular, then we have
\begin{align*}
D(\mathfrak{h}^-) &= \{ \frac{1}{2}(\gamma_{i+1}-\gamma_{i}) \ | \ 1 \le i \le r-1 \} \cup \{0\}\\
\Delta_I \cap \Delta^+ (\mathfrak{h}^-) &= \{ \frac{1}{2}(\gamma_{j}-\gamma_{i}) \ | \ 1 \le i \le j \le r \}\\
\Delta^+ \setminus \Delta_I (\mathfrak{h}^-) &= \{ \frac{1}{2}(\gamma_{j}+\gamma_{i}) \ | \ 1 \le i \le j \le r \}.
\end{align*}
	\item If $(L_I, \mathfrak{n}_I^{+})$ is not regular, then we have
\begin{align*}
D(\mathfrak{h}^-) &= \{ \frac{1}{2}(\gamma_{i+1}-\gamma_{i}) \ | \ 1 \le i \le r-1 \} \cup \{-\frac{1}{2} \gamma_r\} \cup \{0\}\\
\Delta_I \cap \Delta^+ (\mathfrak{h}^-) &= \{ \frac{1}{2}(\gamma_{j}-\gamma_{i}) \ | \ 1 \le i \le j \le r \} \cup \{ -\frac{1}{2} \gamma_i \ | \ 1 \le i \le r \}\\
\Delta^+ \setminus \Delta_I (\mathfrak{h}^-) &= \{ \frac{1}{2}(\gamma_{j}+\gamma_{i}) \ | \ 1 \le i \le j \le r \} \cup \{ \frac{1}{2} \gamma_i \ | \ 1 \le i \le r \}.\\
\end{align*}
\end{enumerate}
\end{lem}

\begin{lem}\label{lem uniqueness}
If $(L_I, \mathfrak{n}_I^+)$ is not regular, then $j \in I$ such that
$\alpha_j |_{\mathfrak{h}^-} = -\frac{\gamma_r}{2} |_{\mathfrak{h}^-}$ is unique.
\end{lem}

\begin{proof}
Assume that for $j_1 \neq j_2$ we have $\alpha_{j_1} = \alpha_{j_2} = -\frac{\gamma_r}{2}$ on $\mathfrak{h}^-$.
Let $h \in \mathfrak{h}^-$.
Since $(\alpha_{j_1} + \alpha_{j_2})(h) = -\gamma_r(h)$, $\alpha_{j_1} + \alpha_{j_2} \notin \Delta$ by Lemma \ref{lem of gamma}.
So we have $(\alpha_{j_1}, \alpha_{j_2}) = 0$.
Since $(\gamma_{r}, \alpha_{j_1})>0$ and $(\gamma_r+\alpha_{j_1}, \alpha_{j_2}) >0$, $\beta = \gamma_r + \alpha_{j_1} + \alpha_{j_2} \in \Delta$.
In particular $\beta \in \Delta^+ \setminus \Delta_I$.
For any $1 \le k \le r$, $(\beta - \gamma_k)(h) = -\gamma_k(h)$, hence
$\beta - \gamma_k \in \sum_{i \in I} \mathbb{Z}_{\ge 0} \alpha_{i}$ is not a root. Clearly $\beta + \gamma_k \notin \Delta$.
Therefore $\beta \in \Gamma_{r+1}$.
But by definition $\Gamma_{r+1} = \emptyset$.
Thus $j_1 = j_2$.
\end{proof}

\begin{lem}\label{lem-w I lambda r}
$w_I \lambda_r = w_0 \varpi_{i_0} - \varpi_{i_0}$.
\end{lem}

\begin{proof}
Assume that $(L_I, \mathfrak{n}_I^+)$ is regular.
Then we have $w_I(w_0 \varpi_{i_0} - \varpi_{i_0}) = -2 w_I \varpi_{i_0} = -2 \varpi_{i_0}$.
By using Lemma \ref{lem of gamma} we can check $\lambda_r(h_i) = -2 \delta_{i, i_0}$ easily (see Wachi \cite{Wachi}).
Hence $\lambda_r = -2 \varpi_{i_0}$, and the statement holds.
Next we assume that $(L_I, \mathfrak{n}_I^+)$ is not regular.
Then there exists an index $j_0 \in I$ such that $\alpha_{j_0} = -\frac{\gamma_r}{2}$ as functions on $\mathfrak{h}^-$.
By Lemma \ref{lem uniqueness} this index $j_0$ is unique.
Similarly to the regular case we have $\lambda_r = \varpi_{j_0} - 2 \varpi_{i_0}$.
Let $j'_0 \in I$ such that $w_I w_0 \alpha_{i_0} = \alpha_{j'_0}$.
Then we have $w_I(w_0 \varpi_{i_0} - \varpi_{i_0}) = \varpi_{j'_0} - 2 \varpi_{i_0}$, and 
we can check $\gamma_r + \alpha_{j'_0} \in \Delta$ by the direct calculation.
Hence $\alpha_{j'_0} = -\frac{\gamma_r}{2}$ on $\mathfrak{h}^-$ from Lemma \ref{lem of gamma}. Therefore we have $j'_0 = j_0$.
\end{proof}

The following fact is known (see \cite{Johnson}, \cite{Schmid}, \cite{Takeuchi}).

\begin{lem}\label{lem-decomp. of Un- in l_I-mod.}
As an $\mathrm{ad}(\mathfrak{l}_I)$-module, $U(\mathfrak{n}_I^{-})$ is multiplicity free, and
\[
U(\mathfrak{n}_I^-) = \bigoplus_{\mu \in \sum_{i=1}^{r} \mathbb{Z}_{\ge 0} \lambda_{i}} I(\mu),
\]
where $I(\mu)$ is an irreducible $\mathfrak{l}_I$-submodule of $U(\mathfrak{n}_I^-)$ with highest weight $\mu$.
\end{lem}
Let $f_{i} \in U(\mathfrak{n}_{I}^-)$ be the highest weight vector of $I(\lambda_i)$.
Since $U(\mathfrak{n}_I^{-})$ is naturally identified with the symmetric algebra $S(\mathfrak{n}_I^{-})$, we can determine the degree of $f \in U(\mathfrak{n}_{I}^-)$.
If $f \in U(\mathfrak{n}_I^{-})$ is a weight vector with weight $\mu \in - d \alpha_{i_0} + \sum_{i \in I} \mathbb{Z}_{\le 0} \alpha_i$, then $f$ is homogeneous and $\deg f = d$.
In particular $\deg f_i = i$.

For $\mu \in \mathfrak{h}_{I, +}^*$ we take the lowest weight vector $v_{w_{0} \mu}$ of $V(\mu)$.
Then the $U_{R_I}(\mathfrak{g})$-module $U_{R_I}(\mathfrak{g}) \otimes_{U_{R_I}(\mathfrak{p}_I^+)} (R_{I, c} \otimes_{\mathbb{C}} V(\mu))$ is generated by
$1 \otimes 1_{c} \otimes v_{w_{0} \mu}$.
There exists a non-zero element $u_{\mu} \in U_{R_I}(\mathfrak{n}_I^{-})$
such that $\psi_{\mu}(1 \otimes 1_{c} \otimes v_{w_{0} \mu}) =
u_{\mu} \otimes 1_{c + \mu}$, where $\psi_{\mu}$ is a $U_{R_I}(\mathfrak{g})$-module homomorphism defined in Section \ref{sc-UGV}.
Since $y (1 \otimes 1_c \otimes v_{w_{0} \mu}) = 0$ for any $y \in \mathfrak{l}_I \cap \mathfrak{n}^-$,
$u_{\mu} \in U_{R_I}(\mathfrak{n}_I^{-})$ is a lowest weight vector with weight $w_{0} \mu - \mu$ as an $\mathrm{ad}(\mathfrak{l}_I)$-module.
By Lemma \ref{lem-decomp. of Un- in l_I-mod.}
such a lowest weight vector is unique up to constant multiple.
Therefore $u_{\mu} = a_{\mu} u_{\mu}^0$ where $a_{\mu} \in R_I \setminus \{0\}$
and $u_{\mu}^0 \in U(\mathfrak{n}_I^-)$ is the unique lowest weight vector with
weight $w_{0} \mu - \mu$.
 
If $x (1 \otimes 1_{c} \otimes v_{w_{0} \mu}) = 0$ for $x \in U_{R_I}(\mathfrak{g})$, then we have $x u_{\mu}^0 \otimes 1_{c + \mu} = 0$ since $a_{\mu} \neq 0$.
Hence we can define a $U_{R_I}(\mathfrak{g})$-module homomorphism $\psi_{\mu}^{0}$ from $U_{R_I}(\mathfrak{g}) \otimes_{U_{R_I}(\mathfrak{p}_I^+)} (R_{I, c} \otimes_{\mathbb{C}} V(\mu))$ to $M_{R_I}(c + \mu)$ 
by 
\[
\psi_{\mu}^{0}(x (1 \otimes 1_{c} \otimes v_{w_{0} \mu})) =
x u_{\mu}^0 \otimes 1_{c + \mu}
\]
for any $x \in U_{R_I}(\mathfrak{g})$.
We set $\xi_{\mu}^0 = \psi_{\mu}^0 \iota_{\mu} \in \Xi_{\mu}$.

Conversely, from the uniqueness of $u_{\mu}^0$
we have
\[
\psi(1 \otimes 1_c \otimes v_{w_0 \mu}) = a u_{\mu}^0 \otimes 1_{c + \mu} = a \psi_{\mu}^0 (1 \otimes 1_c \otimes v_{w_0 \mu}) \quad (a \in R_I)
\]
for any $\psi \in
\mathrm{Hom}_{U_{R_I}(\mathfrak{g})}\big(U_{R_I}(\mathfrak{g}) \otimes_{U_{R_I}(\mathfrak{p}_I^+)} (R_{I, c} \otimes_{\mathbb{C}} V(\mu)), \ M_{R_I}(c + \mu)\big)$.
Therefore we have the following.

\begin{prop}\label{prop-commu diag of comm parabolic}
Let $\mu \in \mathfrak{h}_{I, +}^*$.
We have $\Xi_{\mu} = R_I \xi_{\mu}^0$.
\end{prop}
We call the above homomorphism $\psi_{\mu}^0$ the minimal map in this paper.

Let $\tilde{f}_r$ be the lowest weight vector of the irreducible $\mathfrak{l}_I$-module generated by $f_r$.
\begin{prop}\label{prop-two differ.op.}
Let $\mu = m \varpi_{i_0} \in \mathfrak{h}_{I, +}^*$.
Under the identification $\exp : \mathfrak{n}_I^- \simeq N_I^-$ we have
\[
(\tilde{P}_{\mu}(\psi_{\mu}^0) \varphi)|_{\mathfrak{n}_I^-} = \tilde{f}_{r}(\partial)^m (\varphi|_{\mathfrak{n}_I^-}).
\]
\end{prop}

\begin{proof}

Let $\{ v_{i} \}_{0 \le i \le n}$ be a basis of $V(\mu)$ consisting of weight vectors such that $v_{n}$ has the lowest weight $w_0 \mu$.
We denote the dual basis by $\{v_{i}^{*}\}$.
We define elements $Y'_i \in U_{R_I}(\mathfrak{n}_I^-)$ by $\psi_{\mu}^0 (1 \otimes 1_{c} \otimes v_{i}) = Y'_{i} \otimes 1_{c + \mu}$.
Set $Y_i = \pi(Y'_i)$
Then we have 
\begin{align*}
(\tilde{P}_{\mu}(\psi_{\mu}^0)\varphi)(g) =  \sum_{i = 0}^n \langle g v_{i}^{*}, v_{n} \rangle (R(Y_{i})\varphi)(g).
\end{align*}
For $g \in N_I^{-}$ we have $\langle g v_{i}^{*}, v_{n} \rangle = \delta_{i,n}$. Therefore it is sufficient to show that
\begin{align}\label{eq-diff.op.}
R(Y_n) = \tilde{f}_{r}^m (\partial)
\end{align}
By the definition of $\psi_{\mu}^0$, $Y_{n}$ is the lowest weight vector of $\mathrm{ad}(\mathfrak{l}_I)$-module $U(\mathfrak{n}_I^{-})$ with weight $w_{0} \mu - \mu = m (w_0 \varpi_{i_0} - \varpi_{i_0})$.
By Lemma \ref{lem-w I lambda r} the weight of $\tilde{f_r}$ is $w_0 \varpi_{i_0} - \varpi_{i_0}$.
Hence we have $Y_{n} = \tilde{f}_{r}^m$ up to constant multiple.
Since $\mathfrak{n}_I^-$ is commutative, we have $R(y) = y(\partial)$ for any $y \in U(\mathfrak{n}_I^-)$.
Hence the equation \eqref{eq-diff.op.} holds.
\end{proof}

For $1 \le p \le r = r(\mathfrak{g})$ we set
\begin{align*}
\Delta_{(p)}^+ = \{\beta \in \Delta^+ \ | \ \beta|_{\mathfrak{h}^-} = \frac{\gamma_j + \gamma_k}{2}|_{\mathfrak{h}^-} \ \textrm{ for some } 1 \le j \le k \le p \}.
\end{align*}
By Lemma \ref{lem of gamma} we have $\Delta_{(p)}^+ \subset \Delta^+ \setminus \Delta_I$.
We define subspaces $\mathfrak{n}_{(p)}^{\pm}$ of $\mathfrak{g}$ by $\mathfrak{n}_{(p)}^{\pm} = \sum_{\beta \in \Delta_{(p)}^+} \mathfrak{g}_{\pm \beta}$.
Set $\mathfrak{l}_{(p)} =  [\mathfrak{n}_{(p)}^{+}, \mathfrak{n}_{(p)}^{-}]$ and $I_{(p)} = \{i \in I \ | \ \mathfrak{g}_{\alpha_{i}} \subset \mathfrak{l}_{(p)}\}$.
Then we have the following.
\begin{lem}[{\rm see Wallach \cite{Wallach} and Wachi \cite{Wachi}}]\label{lem-sub PV}
Set $\mathfrak{g}_{(p)} = \mathfrak{n}_{(p)}^{-} \oplus \mathfrak{l}_{(p)} \oplus \mathfrak{n}_{(p)}^{+}$. Then $\mathfrak{g}_{(p)}$ is a simple subalgebra of $\mathfrak{g}$ with simple roots $\{\alpha_{i_0}\} \sqcup \{\alpha_i \ | \ i \in I_{(p)}\}$, and the pair $(\mathfrak{g}_{(p)}, i_0)$ is of regular commutative parabolic type.
For any $1 \le j \le p$ we have $f_j \in U(\mathfrak{n}_{(p)}^-)$,
and $f_{p}$ is a basic relative invariant of $(L_{(p)}, \mathfrak{n}_{(p)}^+)$,
where $L_{(p)}$ is the subgroup of $L_I$ corresponding to $\mathfrak{l}_{(p)}$.
\end{lem}
This fact will be used in the subsequent sections.

\section{Regular Type}
In this section we assume that the prehomogeneous vector space $(L_I, \mathfrak{n}_I^{\pm})$ is regular.
We take $\gamma_{i}, \lambda_{i}$ and $f_{i}$ ($1 \le i \le r = r(\mathfrak{g})$) as in Section
\ref{sc-commuta}.
Then we have $w_{0} \varpi_{i_0} = - \varpi_{i_0}$ and the highest weight vector $f_{r} \in U(\mathfrak{n}_I^-) \simeq \mathbb{C}[\mathfrak{n}_I^+]$ is the unique basic relative invariant of $(L_I, \mathfrak{n}_I^+)$ with character $2 \varpi_{i_0}$.
In particular $f_{r} \in U(\mathfrak{n}_I^-)$ is also the lowest weight vector as an $\mathfrak{l}_I$-module.

\begin{prop}\label{prop- xi and b}
Let $b(s)$ be the $b$-function of the basic relative invariant of $(L_I, \mathfrak{n}_I^-)$.
Then for $m \in \mathbb{Z}_{> 0}$ we have
\[
\xi_{m \varpi_{i_0}}^0 (s \varpi_{i_0}) = b(s + m -1) b(s + m -2) \cdots b(s)
\] up to constant multiple.
\end{prop}

\begin{proof}
For any $l \in L_I$ and $n \in \mathfrak{n}_I^-$ we have
\begin{align*}
\tilde{f}^{\varpi_{i_0}}(l \exp(n) l^{-1}) = (w_0 \varpi_{i_0}-\varpi_{i_0})(l)\tilde{f}^{\varpi_{i_0}}(\exp(n)) = -2 \varpi_{i_0}(l)\tilde{f}^{\varpi_{i_0}}(\exp(n)).
\end{align*}
Thus $\tilde{f}^{\varpi_{i_0}} |_{\mathfrak{n}_I^-}$ is the basic
relative invariant of $(L_I, \mathfrak{n}_I^-)$ under the identification
$\mathfrak{n}_I^- \simeq N_I^-$.
Hence we have
\[
f_r(\partial)^m \tilde{f}^{(m + s)\varpi_{i_0}} |_{\mathfrak{n}_I^-}
= f_r(\partial)^m (\tilde{f}^{\varpi_{i_0}} |_{\mathfrak{n}_I^-})^{s + m}
= b(s + m -1) b(s + m -2) \cdots b(s) \tilde{f}^{s \varpi_{i_0}} |_{\mathfrak{n}_I^-}.
\]
From Corollary \ref{cor-b-function} 
we have
\[
\tilde{P}_{m \varpi_{i_0}}(\psi_{m \varpi_{i_0}}^0)
\tilde{f}^{(s+m) \varpi_{i_0}} = \xi_{m \varpi_{i_0}}^{0}(s \varpi_{i_0}) \tilde{f}^{s \varpi_{i_0}}.
\]
Therefore the statement holds
by Proposition \ref{prop-two differ.op.}.
\end{proof}

In the rest of this section we shall show that $\xi_{\varpi_{i_0}} = \xi_{\varpi_{i_0}}^0$ up to constant multiple.

\begin{lem}
For any $1 \le j \le r$ we have $w_I \gamma_{j} = \gamma_{r-j+1}$.
\end{lem}

\begin{proof}
By Lemmas \ref{lem of gamma-0} and \ref{lem of gamma} we have $\gamma_r + \alpha_i \notin \Delta$ for any $i \in I$.
Therefore $\gamma_r$ is the highest weight of the irreducible $\mathrm{ad}(\mathfrak{l}_I)$-module $\mathfrak{n}_I^+$.
Since $\alpha_{i_0} = \gamma_1$ is the lowest weight of $\mathfrak{n}_I^+$,
$w_I \gamma_1 = \gamma_r$.
Let $1 < i \le [\frac{r}{2}]$.
Assume that $w_I \gamma_j = \gamma_{r-j+1}$ for $1 \le j \le i-1$.
Since $j \neq i$, we have $\gamma_{r-i+1} \pm w_I \gamma_j = \gamma_{r-i+1} \pm \gamma_{r-j+1} \notin \Delta \cup \{0\}$.
Hence $w_I \gamma_{r-i+1} \pm \gamma_j \notin \Delta \cup \{0\}$, and we have $w_I \gamma_{r-i+1} \in \Gamma_i$.
By definition $\gamma_i \le w_I \gamma_{r-i+1}$.
So we have $w_I \gamma_i \ge \gamma_{r-i+1}$.
Let us show that $w_I \gamma_i \le \gamma_{r-i+1}$.
By Lemma \ref{lem of gamma} there exist $\gamma_k$ and $\gamma_l$ such that $k \le l$ and $w_I \gamma_i(h) = \frac{\gamma_k + \gamma_l}{2}(h)$ for any $h \in \mathfrak{h}^-$.
In particular $(w_I \gamma_i, \gamma_l)>0$.
Now we have $(w_I \gamma_i, \gamma_m) = (\gamma_i, w_I \gamma_m) = (\gamma_i, \gamma_{r-m+1}) = 0$ for $r-i+2 \le m \le r$.
Hence $l \le r-i+1$.
Since $(w_I \gamma_i, \gamma_l)>0$, $\gamma_l - w_I \gamma_i \in \Delta \cup \{0\}$.
For $h \in \mathfrak{h}^-$ we have $(\gamma_l - w_I \gamma_i)(h)=\frac{\gamma_l-\gamma_k}{2}(h)$.
By Lemma \ref{lem of gamma} $\gamma_l - w_I \gamma_i \in \Delta^+ \cup \{0\}$.
Therefore we have $w_I \gamma_i \le \gamma_l \le \gamma_{r-i+1}$.
\end{proof}

Hence the lowest weight $w_I \lambda_{r-1}$ of the irreducible component
$I(\lambda_{r-1})$ of $U(\mathfrak{n}_I^{-})$ is $\lambda_{r} + \alpha_{i_0}$.

\begin{lem}\label{lem-act on f_p}
For any $1 \le p \le r = r(\mathfrak{g})$ we have
\[
e_{i_0} f_{p} \otimes 1_{c + \mu} \in U_{R_I}(\mathfrak{l}_I \cap \mathfrak{n}^-) (f_{p-1} \otimes 1_{c + \mu})
\subset M_{R_I}(c + \mu),
\]
where $e_{i_0} \in \mathfrak{g}_{\alpha_{i_0}} \setminus \{0\}$.
\end{lem}

\begin{proof}
By Lemma \ref{lem-sub PV} it is sufficient to show that the statement 
holds for $p = r$.
We define $y \in U_{R_I}(\mathfrak{n}_I^-)$ by
\[
e_{i_0} (f_r \otimes 1_{c + \mu}) = y \otimes 1_{c + \mu}.
\]
Since $f_{r}$ is the lowest weight vector of the $\mathrm{ad}(\mathfrak{l}_I)$-module $U(\mathfrak{n}_I^-)$ and $[e_{i_0}, \mathfrak{l}_I \cap \mathfrak{n}^-] = \{0\}$, $y$ is the lowest weight vector as an $\mathrm{ad}(\mathfrak{l}_I)$-module.
Moreover the weight of $y$ is $\lambda_{r} + \alpha_{i_0} = w_I \lambda_{r-1}$,
which is the lowest weight of the irreducible component $I(\lambda_{r-1}) = 
\mathrm{ad}(U(\mathfrak{l}_I))f_{r-1}$.
Therefore $y \otimes 1_{c + \mu} \in U_{R_I}(\mathfrak{l}_I \cap \mathfrak{n}^-) (f_{r-1} \otimes 1_{c + \mu})$.
\end{proof}

\begin{cor}\label{cor-f_r}
Let $u \in U(\mathfrak{n}^+)$ with weight $k \alpha_{i_0} + \sum_{i \in I} m_i \alpha_i$.
Then we have
\[
u f_r \otimes 1_{c + \mu} \in U_{R_I}(\mathfrak{l}_I \cap \mathfrak{n}^-)
(f_{r-k} \otimes 1_{c + \mu}).
\]
\end{cor}
\begin{proof}
We shall show the statement by the induction on $k$.
If $k = 0$, then the statement is clear.
Assume that $k > 0$, and the statement holds for $k-1$.
We write $u = \sum_{j} u_j e_{i_0} u'_j$, where $u_j \in U(\mathfrak{l}_I \cap \mathfrak{n}^+)$ and $u'_j \in U(\mathfrak{n}^+)$.
Then the weight of $u'_j$ is in $(k-1) \alpha_{i_0} + \sum_{i \in I} \mathbb{Z}_{\ge 0} \alpha_i$, and hence we have
\[
u f_r \otimes 1_{c + \mu} \in \sum_{j} u_j e_{i_0} U_{R_I}(\mathfrak{l}_I \cap \mathfrak{n}^-) (f_{r-k+1} \otimes 1_{c + \mu}) \subset U_{R_I}(\mathfrak{l}_I)(e_{i_0} f_{r-k+1} \otimes 1_{c + \mu}).
\]
Here note that $[e_{i_0}, U_{R_I}(\mathfrak{l}_I \cap \mathfrak{n}^-)] = 0$.
By Lemma \ref{lem-act on f_p} we have 
\[
e_{i_0} f_{r-k+1} \otimes 1_{c + \mu}
\in U_{R_I}(\mathfrak{l}_I \cap \mathfrak{n}^-) (f_{r-k} \otimes 1_{c + \mu}).
\]
Therefore we obtain
\[
u f_r \otimes 1_{c + \mu} \in U_{R_I}(\mathfrak{l}_I)(f_{r-k} \otimes 1_{c + \mu}) = U_{R_I}(\mathfrak{l}_I \cap \mathfrak{n}^-)(f_{r-k} \otimes 1_{c + \mu}).
\]
\end{proof}

\begin{thm}\label{thm-explicit of xi^0}
We have 
$\xi_{\varpi_{i_0}} =\prod_{j=1}^r p_{\varpi_{i_0}, \lambda_j + \varpi_{i_0}} \in
\mathbb{C}^{\times} \xi_{\varpi_{i_0}}^0$, where $\mathbb{C}^{\times} = \mathbb{C} \setminus \{0\}$.
\end{thm}

\begin{proof}
Let $v_{- \varpi_{i_0}}$ be the lowest weight vector of $V(\varpi_{i_0})$.
Since $f_r$ is the lowest weight vector of $U(\mathfrak{n}_I^-)$ with weight $-2 \varpi_0$, we have $\psi_{\varpi_{i_0}}^{0}(1 \otimes 1_{c} \otimes v_{-\varpi_{i_0}}) = f_r \otimes 1_{c+\varpi_{i_0}}$.
It is clear that
\[
\varpi_{i_0} - w_0 \varpi_{i_0} = 2 \varpi_{i_0} = - \lambda_{r} \in r \alpha_{i_0} + \sum_{i \in I} \mathbb{Z}_{\ge 0} \alpha_{i}.
\]
Set $P(j) = \{ \lambda \ | \ \varpi_{i_0} - \lambda \in j \alpha_{i_0} + \sum_{i \in I} \mathbb{Z}_{\ge 0} \alpha_{i} \}$.
We define an $\mathfrak{l}_I$-submodule $V_j$ of $V(\varpi_{i_0})$
by
\[
V_j = \bigoplus_{\lambda \in P(j)} V(\varpi_{i_0})_{\lambda}
\]
(cf. Section \ref{sc-UGV}).
Note that $V_j \neq 0$ for $0 \le j \le r$.
We take the irreducible decomposition of $V_j$
\[
V_j = \tilde{W}(\eta_{j,1})
\oplus \cdots \oplus \tilde{W}(\eta_{j,t_j}),
\]
where $\tilde{W}(\eta)$ is an irreducible $\mathfrak{l}_I$-module with
highest weight $\eta$.
Let $v_{j,k}$ be the highest weight vector of $\tilde{W}(\eta_{j,k})$. 
There exists an element $u_{j,k} \in U(\mathfrak{n}^+)$ such that $u_{j,k} v_{- \varpi_{i_0}} = v_{j,k}$.
Then the weight of $u_{j,k}$ is in $(r-j) \alpha_{i_0} + \sum_{i \in I} \mathbb{Z}_{\ge 0} \alpha_{i}$.
By Corollary \ref{cor-f_r} we have
\[
\psi_{\varpi_{i_0}}^{0}(1 \otimes 1_c \otimes v_{j,k}) =
u_{j,k} \psi_{\varpi_{i_0}}^{0}(1 \otimes 1_c \otimes v_{- \varpi_{i_0}}) =
u_{j,k} f_r \otimes 1_{c + \varpi_{i_0}} \in
U_{R_I}(\mathfrak{l}_I \cap \mathfrak{n}^-) f_{j} \otimes 1_{c + \varpi_{i_0}}.
\]
Since $v_{j,k}$ is the highest weight vector, we have
\[
\psi_{\varpi_{i_0}}^{0}(1 \otimes 1_c \otimes v_{j,k}) \in R_I( f_{j} \otimes 1_{c + \varpi_{i_0}}).
\]
In particular $\eta_{j,k} = \lambda_{j} + \varpi_{i_0}$, and
the irreducible decomposition of $V(\varpi_{i_0})$ as an $\mathfrak{l}_I$-module is given by
\[
V(\varpi_{i_0}) = \bigoplus_{j = 0}^r \tilde{W}(\lambda_{j} + \varpi_{i_0})^{\oplus N_{j}},
\]
where we set $\lambda_0 = 0$.
Therefore we have $\xi_{\varpi_{i_0}} = \prod_{j=1}^r p_{\varpi_{i_0}, \lambda_j + \varpi_{i_0}}$,
which is regarded as a polynomial function on $\mathbb{C} \varpi_{i_0}$.
Since $\deg p_{\varpi_{i_0}, \lambda_j + \varpi_{i_0}} = 1$ for $j \ge 1$,
we have $\deg \xi_{\varpi_{i_0}} = r$.
Now we have $\deg \xi_{\varpi_{i_0}}^0 = \deg b(s) = \deg f_r = r$.
From Proposition \ref{prop-commu diag of comm parabolic}
we have $\xi_{\varpi_{i_0}} \in \mathbb{C}^{\times} \xi_{\varpi_{i_0}}^0$, 
hence the statement holds.
\end{proof}

For $1 \le i < j \le r$ we set $c_{i,j} = \sharp \{ \alpha \in \Delta_I \cap \Delta^+ \ | \ \alpha |_{\mathfrak{h}^-} = \frac{\gamma_j - \gamma_i}{2}|_{\mathfrak{h}^-}\}$.
It is known that $c_{i,j} = \sharp \{ \alpha \in \Delta^+ \setminus \Delta_I \ | \ \alpha |_{\mathfrak{h}^-} = \frac{\gamma_j + \gamma_i}{2}|_{\mathfrak{h}^-}\}$ and this number is independent of $i$ or $j$ (see \cite{Wallach2}).
Set $c_0 = c_{i,j}$.
Then we have $(2 \rho, \gamma_j) = d_0 (1 + c_0 (j-1))$, where $d_0 = (\alpha_{i_0}, \alpha_{i_0})$.
In particular $(2 \rho, \lambda_j)= -j d_0 (1 + \frac{j-1}{2}c_0)$.
Since $(\gamma_i, \gamma_j) = \delta_{i,j} d_0$, we have $(\varpi_{i_0}, \varpi_{i_0}) = (\lambda_j+\varpi_{i_0}, \lambda_j+\varpi_{i_0})$ for $1 \le j \le r$.
Hence we have
\[
p_{\varpi_{i_0}, \lambda_j+\varpi_{i_0}}(s \varpi_{i_0}) = -2 (s \varpi_{i_0} + \rho, \lambda_{j}) = j d_0 (s + 1 + \frac{j-1}{2}c_0).
\]

\section{Non-regular Type}

Assume that the prehomogeneous vector space $(L_I, \mathfrak{n}_I^+)$ is not regular. 
We take $\gamma_i$, $\lambda_i$ and $f_i$ $(1 \le i \le r = r(\mathfrak{g}))$ as in Section \ref{sc-commuta}.
For $\mu = m \varpi_{i_0} \in \mathfrak{h}_{I, +}^*$ we denote by $\tilde{v}_{\mu}$ the highest weight vector of the irreducible $\mathfrak{l}_I$-submodule of $V(\mu)$ generated by the lowest weight vector of $V(\mu)$.
The weight of $\tilde{v_{\mu}}$ is $w_I \mu$.
We take $u \in U_{R_I}(\mathfrak{n}_I^-)$ as $\psi_{\mu}^0(1 \otimes 1_c \otimes \tilde{v}_{\mu}) = u \otimes 1_{c+\mu}$.
By definition of $\psi_{\mu}^0$ we have $u \in U(\mathfrak{n}_I^-)$.
Moreover $u$ is the highest weight vector of $U(\mathfrak{n}_I^-)$ with weight $w_I w_0 \mu - \mu = w_I(w_0 \mu - \mu)$.
By Lemma \ref{lem-w I lambda r} we have $w_I(w_0 \mu - \mu) = m \lambda_r$.
Therefore we have $u = f_{r}^m$.
Set $\xi_{\mu}^0 = \psi_{\mu}^0 \iota_{\mu} \in R_I$.

We define subalgebras $\mathfrak{g}_{(r)}$, $\mathfrak{l}_{(r)}$ and $\mathfrak{n}_{(r)}^{\pm}$ of $\mathfrak{g}$ as in Lemma \ref{lem-sub PV}.
We set $\tilde{\mathfrak{p}}^+ = \mathfrak{l}_{(r)} \oplus \mathfrak{n}_{(r)}^+$.
We denote by $\tilde{V}(\mu)$ the irreducible $\mathfrak{g}_{(r)}$-module with highest weight $\mu$.
Let $\tilde{I}_0$ be an index set of simple roots of $\mathfrak{g}_{(r)}$, that is, $\tilde{I}_0 = I_{(r)} \sqcup \{i_0\}$ (see Lemma \ref{lem-sub PV}).
We set $\tilde{I} = I_{(r)}$ and $\tilde{\mathfrak{g}} = \mathfrak{g}_{(r)}$ for simplicity.
Let $\tilde{R}$ be an enveloping algebra of $\sum_{i \in \tilde{I}_0} \mathbb{C}h_i / \sum_{i \in \tilde{I}} \mathbb{C} h_i$.
Since we have the canonical identification $R_I \simeq \tilde{R}$, a $U_{R_{I}}(\tilde{\mathfrak{g}})$-submodule 
\[
\tilde{M}(c+\mu) = U_{R_{I}}(\tilde{\mathfrak{g}}) \otimes_{U_{R_{I}}(\tilde{\mathfrak{p}}^+)} R_{I, c+\mu}
\]
of $M_{R_I}(c+\mu)$ is a generalized universal Verma module associated with $\tilde{\mathfrak{g}}$.
We define an element $\tilde{\xi}_{\mu}^0$ of $\tilde{R} \simeq R_{I}$ by the multiplication map on $\tilde{M}(c+\mu)$ induced by the minimal map
\[
\tilde{\psi}_{\mu}^0 : U_{R_{I}}(\tilde{\mathfrak{g}}) \otimes_{U_{R_{I}}(\tilde{\mathfrak{p}}^+)} (R_{I, c} \otimes_{\mathbb{C}} \tilde{V}(\mu)) \to \tilde{M}(c+\mu).
\]
Then we have the following.
\begin{prop}\label{prop-non-regular}
\begin{enumerate}
	\item Under the identification $\tilde{R} \simeq R_{I}$ we have $\xi_{\mu}^0 = \tilde{\xi}_{\mu}^0$ for $\mu \in \mathfrak{h}_{I, +}^*$. \label{reg and nonreg}
	\item $\xi_{\varpi_{i_0}} \in \mathbb{C}^{\times} {\xi}_{\varpi_{i_0}}^0$.
\end{enumerate}
\end{prop}
\begin{proof}
\begin{enumerate}
	\item We have $U(\tilde{\mathfrak{g}}) {v}_{\mu} \simeq \tilde{V}(\mu)$, and $\tilde{v}_{\mu}$ is its lowest weight vector.
Since we have $\tilde{\psi}_{\mu}^0(1 \otimes 1_{c} \otimes \tilde{v}_{\mu}) = f_r^{m} \otimes 1_{c + \mu}$, the restriction $\psi_{\mu}^0$ on $U_{R_{{I}}}(\tilde{\mathfrak{g}}) \otimes_{U_{R_{{I}}}(\tilde{\mathfrak{p}}^+)} (R_{{I}, c} \otimes_{\mathbb{C}} U(\tilde{\mathfrak{g}}) v_{\mu})$ is $\tilde{\psi}_{\mu}^0$.
Hence we have $\xi_{\mu}^0 \otimes 1_{c+\mu} = \psi_{\mu}^0(1 \otimes 1_{c} \otimes v_{\mu}) = \tilde{\psi}_{\mu}^0(1 \otimes 1_{c} \otimes v_{\mu}) = \tilde{\xi}_{\mu}^0 \otimes 1_{c+\mu}$.
	\item Since the pair $(\tilde{\mathfrak{g}}, i_0)$ is of regular type, we have $\deg \tilde{\xi}_{\varpi_{i_0}}^0 = r$ (see the proof of Theorem \ref{thm-explicit of xi^0}).
Similarly to the proof of Theorem \ref{thm-explicit of xi^0} we can show that $\deg \xi_{\varpi_{i_0}} = r$.
By (\ref{reg and nonreg}) we have $\deg \xi_{\varpi_{i_0}}^0 = \deg {\xi}_{\varpi_{i_0}}$.
Since $\Xi_{\varpi_{i_0}} = R_I {\xi}_{\varpi_{i_0}}^0$ and ${\xi}_{\varpi_{i_0}} \in \Xi_{\varpi_{i_0}}$,
we have $\xi_{\varpi_{i_0}} \in \mathbb{C}^{\times} {\xi}_{\varpi_{i_0}}^0$.
\end{enumerate}
\end{proof}

As a result, we have the following.
\begin{thm}\label{generator of Xi_{varpi_{i_0}}}
For any pair $(\mathfrak{g}, i_0)$ of commutative parabolic type , the ideal $\Xi_{\varpi_{i_0}}$ is generated by $\xi_{\varpi_{i_0}}$.
\end{thm}

\section{Irreducibility of Verma Modules}
Let $(L_I, \mathfrak{n}_I^-)$ be a prehomogeneous vector space
of commutative parabolic type. Set $\{ i_0 \} = I_0 \setminus I$.
In this section we give a new proof of the following well-known fact (Suga \cite{Suga}, Gyoja \cite{Gyoja}, Wachi \cite{Wachi}).
\begin{thm}\label{thm-irreducibility}
Let $\lambda = s_{0} \varpi_{i_0} \in \mathfrak{h}_I^*$.
$M_{I}(\lambda)$ is irreducible if and only if
$\xi_{\varpi_{i_0}}^0(\lambda - m \varpi_{i_0}) \neq 0$ for any $m \in \mathbb{Z}_{>0}$.
\end{thm}

We take $f_i \in U(\mathfrak{n}_I^-)$ ($1 \le i \le r = r(\mathfrak{g}$)) as in Section \ref{sc-commuta}.

\begin{lem}\label{lem-decomp.of f_r}
Let $e_{i}^-$ be a nonzero element of $\mathfrak{g}_{-\alpha_i}$.
There exist $y_0, y_1, \ldots, y_t \in \mathfrak{n}_{(r)}^-$ and $ i_1,\ldots, i_t \in I_{(r)}$ such that
\begin{align}\label{eq-decomp f}
f_{r} = \sum_{k=0}^t y_k \mathrm{ad}(e^-_{i_k} \cdots e^-_{i_1}) f_{r-1},
\end{align}
and $\mathrm{ad}(e^-_{i_k}) f_{r-1} = 0$ for $k \ge 2$.
\end{lem}
This lemma is proved by direct calculations for each case. 
In \cite{Kamita} there are explicit decompositions of quantum counterparts $f_{q,r}$ of $f_r$ satisfying the properties of Lemma \ref{lem-decomp.of f_r}.
We can get the decomposition \eqref{eq-decomp f} from the quantum counterpart via $q = 1$.
For example, in the case of type $A$, 
$f_r$ is a determinant and the decomposition \eqref{eq-decomp f} corresponds to a cofactor decomposition.

\begin{prop}\label{prop-f_i->f_{i+1}}
Let $a_{i} \in \mathbb{Z}_{>0}$ and let $a_{i+1}, \ldots, a_{r} \in \mathbb{Z}_{\ge 0}$.
There exists $u \in U(\mathfrak{n}^-)$ such that
\[
u f_i^{a_i} f_{i+1}^{a_{i+1}} \cdots f_r^{a_r} \otimes 1_{\lambda} = f_i^{a_{i}-1} f_{i+1}^{a_{i+1} + 1} \cdots f_{r}^{a_r} \otimes 1_{\lambda}.
\]
\end{prop}

\begin{proof}
By Lemmas \ref{lem-sub PV} and \ref{lem-decomp.of f_r}
we have
\[
f_{i+1} = \sum_{k = 0}^t y_k \mathrm{ad}(u_k) f_{i},
\]
where $y_k \in \mathfrak{n}_{(i+1)}^-$ and
$u_k = e_{j_k}^- \cdots e_{j_1}^-$ such that $j_1, \ldots, j_{t} \in I_{(i+1)}$
and $\mathrm{ad}(e_{j_l}^-) f_{i} = 0$ for $l >1$.
Note that for $p > i$ we have $y_k \in \mathfrak{n}_{(p)}^-$ and
$j_k \in I_{(p)}$.
Since $f_{p}$ is the lowest weight vector of an $\mathrm{ad}(\mathfrak{l}_{(p)})$-module $U(\mathfrak{n}_{(p)}^-)$, we have
\[
\mathrm{ad}(u_k)(f_{i}^{a_i } f_{i+1}^{a_{i+1}} \cdots f_r^{a_r})
= (\mathrm{ad}(u_k) f_{i}^{a_i }) f_{i+1}^{a_{i+1}} \cdots f_r^{a_r}.
\]
If $k \ge 1$, then we have
\begin{align*}
\mathrm{ad}(u_k)(f_{i}^{a_i }) &= a_i \ \mathrm{ad}(e_{j_k}^- \cdots e_{j_{2}}^-)
((\mathrm{ad}(e_{j_1}^-)f_{i}) f_{i}^{a_i - 1})\\
&=  a_i \ (\mathrm{ad}(u_k)f_i)f_{i}^{a_i - 1}.
\end{align*}
Hence for $u = y_0 + a_i^{-1} \sum_{k = 1}^t y_k u_k$,
we have
\begin{align*}
u f_i^{a_i} f_{i+1}^{a_{i+1}} \cdots f_r^{a_r} \otimes 1_{\lambda}
&= y_0 f_i^{a_i} f_{i+1}^{a_{i+1}} \cdots f_r^{a_r} \otimes 1_{\lambda}
+ \sum_{k=1}^t y_k (\mathrm{ad}(u_k)f_i) f_{i}^{a_i - 1} f_{i+1}^{a_{i+1}}
\cdots f_{r}^{a_r} \otimes 1_{\lambda}\\
&= f_i^{a_{i}-1} f_{i+1}^{a_{i+1} + 1} \cdots f_{r}^{a_r} \otimes 1_{\lambda}.
\end{align*}
\end{proof}

\begin{cor}\label{cor-submod}
Let $K (\neq 0)$ be a submodule of $M_I(\lambda)$ for $\lambda \in \mathfrak{h}_I^*$.
We have $f_r^{n} M_I(\lambda) \subset K$ for $n \gg 0$.
\end{cor}

\begin{proof}
If $K = M_I(\lambda)$, then the statement is clear.
Assume that $\{0\} \neq K \subsetneq M_I(\lambda)$.
By Lemma \ref{lem-decomp. of Un- in l_I-mod.} any highest weight vector of $M_I(\lambda)$ as an $\mathfrak{l}_I$-module is given by the following form:
\[
f_1^{a_1} \cdots f_{r}^{a_r} \otimes 1_{\lambda}.
\]
Since $K$ has the highest weight vector as an $\mathfrak{l}_I$-module,
there exists an element $f_1^{a_1} \cdots f_{r}^{a_r} \otimes 1_{\lambda} \in K$
such that $(a_1, \ldots, a_r) \neq 0$.
By Proposition \ref{prop-f_i->f_{i+1}} for $n \gg 0$ there exists
$u \in U(\mathfrak{n}^-)$ such that
\[
f_{r}^{n} \otimes 1_{\lambda} = u (f_1^{a_1} \cdots f_{r}^{a_r} \otimes 1_{\lambda}) \in K.
\]
Hence for any $y \in U(\mathfrak{n}_I^-)$ we have
\[
f_{r}^{n} (y \otimes 1_{\lambda}) = y f_{r}^{n} \otimes 1_{\lambda} \in K,
\]
and the statement holds.
\end{proof}

Let us prove Theorem \ref{thm-irreducibility} by using the commutative diagram
\begin{align}\label{eq-cd-using proof}
\begin{CD}
M_{R_I}(c + \mu) @= M_{R_I}(c + \mu)\\
@V{\iota_{\mu}}VV @VV{\xi_{\mu}^0}V\\
U_{R_I}(\mathfrak{g}) \otimes_{U_{R_I}(\mathfrak{p}_I^+)} (R_{I, c} \otimes_{\mathbb{C}} V(\mu)) @>>\psi_{\mu}^0> M_{R_I}(c + \mu),\\
\end{CD}
\end{align}
where $\mu \in \mathfrak{h}_{I,+}^*$.

Set $\lambda = s_0 \varpi_{i_0}$.
We denote the highest weight vector of $V(\mu)$
by $v_{\mu}$.
Let $\tilde{v}_{\mu}$ be the highest weight vector of the irreducible $\mathfrak{l}_I$-module generated by the lowest weight vector of $\mathfrak{g}$-module $V(\mu)$.
For a positive integer $m$, we set $\mu = \mu_m = m \varpi_{i_0}$.
Considering the functor $\mathbb{C} \otimes_{R_I}(\ \cdot \ )$, where $\mathbb{C}$ has the $R_I$-module structure via $c(h_i) 1 = (\lambda - \mu)(h_i)$,
we obtain the following commutative diagram from \eqref{eq-cd-using proof}:
\begin{align*}
\begin{CD}
M_{I}(\lambda) @= M_{I}(\lambda)\\
@V{\iota_m}VV @VV{\xi_{\mu}^0(\lambda-\mu)}V\\
U(\mathfrak{g}) \otimes_{U(\mathfrak{p}_I^+)} (\mathbb{C}_{I, \lambda-\mu} \otimes_{\mathbb{C}} V(\mu)) @>>{\psi_m^0}> M_{I}(\lambda),\\
\end{CD}
\end{align*}
where $\iota_m(1 \otimes 1_{\lambda}) = 1 \otimes 1_{\lambda-\mu} \otimes v_{\mu}$ and $\psi_m^0(1 \otimes 1_{\lambda-\mu} \otimes \tilde{v}_{\mu}) = f_r^m \otimes 1_{\lambda}$.

Assume that $M_I(\lambda)$ is irreducible.
Since $\psi_m^0 \neq 0$,
we have $\mathrm{Im} \psi_m^0 = M_I(\lambda)$.
The weight space of $U(\mathfrak{g}) \otimes_{U(\mathfrak{p}_I^+)} (\mathbb{C}_{I, \lambda-\mu} \otimes_{\mathbb{C}} V(\mu))$ with weight $\lambda$
is $\mathbb{C}(1 \otimes 1_{\lambda-\mu} \otimes v_{\mu})$, 
hence there exists $a \in \mathbb{C} \setminus \{0\}$ such that
\begin{align*}
1 \otimes 1_{\lambda} = \psi_m^0 (a \otimes 1_{\lambda-\mu} \otimes v_{\mu})
 = a \psi_m^0 \iota_m (1 \otimes 1_{\lambda})
= a \xi_{\mu}^0(\lambda - \mu) \otimes 1_{\lambda} \neq 0.
\end{align*}
By Propositions \ref{prop- xi and b} and \ref{prop-non-regular} we have
\[
\xi_{\mu}^0(\lambda-\mu)
= \xi_{\varpi_{i_0}}^0(\lambda - \varpi_{i_0}) \xi_{\varpi_{i_0}}^0(\lambda - 2 \varpi_{i_0}) \cdots \xi_{\varpi_{i_0}}^0(\lambda - m \varpi_{i_0}).
\]
Therefore we have $\xi_{\varpi_{i_0}}^0(\lambda - m \varpi_{i_0}) \neq 0$ for any $m \in \mathbb{Z}_{>0}$.

Conversely, we assume that $\xi_{\varpi_{i_0}}^0(\lambda - m \varpi_{i_0}) \neq 0$ for any $m \in \mathbb{Z}_{>0}$.
We set
\[
N(m) = U(\mathfrak{g}) \otimes_{U(\mathfrak{p}_I^+)} (\mathbb{C}_{I, \lambda-\mu_m} \otimes_{\mathbb{C}} V(\mu_m)).
\]
Since $\xi_{\mu_m}^0(\lambda-\mu_m)
= \xi_{\varpi_{i_0}}^0(\lambda - \varpi_{i_0}) \xi_{\varpi_{i_0}}^0(\lambda - 2 \varpi_{i_0}) \cdots \xi_{\varpi_{i_0}}^0(\lambda - m \varpi_{i_0}) \neq 0$, we have
\begin{align*}
\psi_{m}^0 (\xi_{\mu_m}^0(\lambda-\mu_m)^{-1} \otimes 1_{\lambda-\mu_m} \otimes v_{\mu_m})
& =  \xi_{\mu_m}^0(\lambda-\mu_m)^{-1} \psi_{m}^0 \iota_{m}(1 \otimes 1_{\lambda})\\
& =  1 \otimes 1_{\lambda}.
\end{align*}
Hence $\psi_{m }^0$ is surjective, and we have an isomorphism
\begin{align*}\label{eq-identify}
N(m)/\mathrm{Ker} \psi_{m }^0
\simeq M_I(\lambda)
: \ \overline{1 \otimes 1_{\lambda-\mu_m} \otimes \tilde{v}_{\mu_m}} \mapsto f_r^m \otimes 1_{\lambda}
\end{align*}
for any $m$.
Under this identification we have
\begin{align*}
\overline{1 \otimes 1_{\lambda - \mu_{n+1}} \otimes \tilde{v}_{\mu_{n+1}}}
= f_{r}^{n+1} \otimes 1_{\lambda} = f_{r}^n \, (f_{r} \otimes 1_{\lambda})
= f_{r}^n \, \overline{1 \otimes 1_{\lambda-\mu_1} \otimes \tilde{v}_{\mu_1}}.
\end{align*}
Let $K \neq 0$ be a submodule of $M_I(\lambda)$.
By Corollary \ref{cor-submod} for $n \gg 0$ we have
\[
\overline{1 \otimes 1_{\lambda-\mu_{n+1}} \otimes \tilde{v}_{\mu_{n+1}}} = 
f_{r}^n \, \overline{1 \otimes 1_{\lambda-\mu_1} \otimes \tilde{v}_{\mu_1}} \in K.
\]
Hence we have 
\begin{align*}
M_I(\lambda)& = N(n+1)/\mathrm{Ker} \psi_{n+1}^0\\
&= U(\mathfrak{g}) \overline{1 \otimes 1_{\lambda-\mu_{n+1}} \otimes \tilde{v}_{\mu_{n+1}}} \subset K.
\end{align*}
Therefore $K = M_I(\lambda)$, and
$M_I(\lambda)$ is irreducible.
We complete the proof of Theorem \ref{thm-irreducibility}.

\end{document}